\documentclass{article}

\usepackage[includeheadfoot,
            bindingoffset=0mm,
            inner   = 25mm,
            top     = 20mm,
            outer   = 25mm,
            bottom  = 10mm,
            paperwidth = 210mm,
            paperheight = 297mm,
            ]{geometry}
\usepackage[margin={2.0cm,0cm},oneside,labelfont={sf,bf},singlelinecheck=false]{caption}

\usepackage{graphicx}

\usepackage[fleqn]{amsmath}
\usepackage{accents}
\usepackage{amssymb}

\usepackage{enumitem}
\usepackage{cite}
\usepackage[perpage]{footmisc}
\usepackage{hyphenat}
\usepackage[english]{babel}
\usepackage[section]{placeins}
\usepackage{color}
\usepackage{upgreek}
\usepackage{listings}

\definecolor{codegreen}{rgb}{0,0.6,0}
\definecolor{codeblue}{rgb}{0,0,0.8}
\definecolor{codegrey}{rgb}{0.5,0.5,0.5}
\usepackage{amsthm}
\theoremstyle{definition}
\newtheorem{theorem}{Theorem}
\newtheorem{remark}{Remark}

\usepackage{titlesec}
\titleformat{\section}[hang]{\Large\bfseries\raggedright\sffamily}{\thesection}{1em}{}
\titleformat{\subsection}[hang]{\large\bfseries\raggedright\sffamily}{\thesubsection}{1em}{}
\titleformat{\subsubsection}[hang]{\normalsize\bfseries\raggedright\sffamily}{\thesubsubsection}{1em}{}
\usepackage{abstract}

\usepackage{authblk}

\newcommand{\transp}{\ensuremath{ ^\mathrm{T} }}

\newcommand{\dif}{\ensuremath{ \mathrm{d} }}

\begin{document}

\title{ \huge\bfseries\sffamily Modelling non-linear control systems using the discrete Urysohn operator }

\author[1]{M. Poluektov}
\author[2]{A. Polar}
\affil[1]{International Institute for Nanocomposites Manufacturing, WMG, University of Warwick, Coventry CV4 7AL, UK}
\affil[2]{Independent software consultant, Duluth, GA, USA}


\maketitle

\setlength{\absleftindent}{2.0cm}
\setlength{\absrightindent}{2.0cm}
\setlength{\absparindent}{0em}
\begin{abstract}
This paper introduces a multiple-input discrete Urysohn operator for modelling non-linear control systems and a technique of its identification by processing the observed input and output signals. It is shown that, due to the nature of the discrete Urysohn operator, the identification problem always has an infinity of solutions, which exactly convert the inputs to the output. The suggested iterative identification procedure, however, leads to a unique solution with the minimum norm, requires only few arithmetic operations with the parameter values and is applicable to a real-time identification, running concurrently with the data reading. The efficiency of the proposed modelling and identification approaches is demonstrated using an example of a non-linear mechanical system, which is represented by a differential equation, and an example of a complex real-world dynamic object.
\end{abstract}

\section{Introduction}
\label{sec:intro}

There is a variety of models used in control system identification \cite{Nelles2001}, e.g. the Volterra series, the Hammerstein, the Wiener-Hammerstein, the Urysohn, the neural network models or the NARMAX model. The Urysohn model is a generalisation of the well-known Hammerstein model\footnote{This is further discussed in appendix \ref{sec:relHammer} of this paper.} and is based on the integral operator of the Urysohn type. These models are the so-called grey-box models and are often differentiated from the physics-based models \cite{Kerschen2006,Ljung2010,Schoukens2019}, as they lack description of the underlying physics of the modelled systems. Nevertheless, these models are often used for a large range of engineering applications.

The major problem with using any generic model of a control object is the identification of the model parameters. There is a number of papers dedicated to different aspects of solving the Urysohn integral equation for a given kernel, e.g. \cite{Balder1984,Angell2010,Alias2016,Curato2016}. However, literature on the identification of the kernel based on known input and output data is very limited. To identify the kernel, the Urysohn model is usually approximated by parallel Hammerstein blocks, e.g. \cite{Gallman1975,Chen1995,Kiselman2005,Harnischmacher2007,Schoukens2011}, by Lagrange polynomials \cite{Makarov1994} or by Stancu polynomials \cite{Makarov2012}. Application-oriented papers, where the Urysohn model identification has been performed, include \cite{Menold1997,Meiler2008,Meiler2009,Moslehpour2015,Moslehpour2016}. 

In contrast to identification of the Urysohn systems, literature on the identification of the Hammerstein systems is vast and most well-known methods include \cite{Bai1998,Bai2004,Liu2007}. However, a detailed discussion of the Hammerstein systems' identification is out of the scope of this paper.

In \cite{Poluektov1990}, it has been suggested to use a discrete Urysohn kernel and to identify it as a grid. Although the idea behind the discrete Urysohn operator is relatively simple, surprisingly, it has not been used in literature since then (to the best knowledge of the authors). This paper picks up the idea of using the discrete Urysohn operator for modelling non-linear control systems and aims at demonstrating that it is an extremely efficient engineering tool --- it describes highly non-linear objects (e.g. with non-monotonous steady-state characteristics) and can be easily identified. This paper builds on \cite{Poluektov1990} by generalising the discrete Urysohn model, which originally has been proposed for a quantised single input only, and by suggesting a more versatile identification technique based on the Kaczmarz iterative method \cite{Kaczmarz1937}. Furthermore, in this paper, the identification problem for the discrete Urysohn operator is considered in detail and the theorem on non-uniqueness of the Urysohn operator (and the structure of the operator with respect to free parameters) is proved.

In \cite{Poluektov1990}, the discrete Urysohn model has been used for modelling the dynamic behaviour of diesel engines and an excellent match with the experimental data has been achieved. However, the identification has been relatively complex, computationally expensive and required the gradual change of the elements of the operator. In this paper, the proposed method identifies the Urysohn kernel based on the observable input-output data only, where the input sequence can be arbitrary, although still covering the desired input range. The proposed method has low computational complexity and can be implemented in a few lines of code.

This paper is organised as follows. In section \ref{sec:theory}, the description of the continuous and the discrete Uryshon models is summarised. Section \ref{sec:ident} holds the major results regarding the solution of the identification problem, with details provided in the appendices. Generalisations of the discrete Urysohn operator and some properties of the Urysohn systems are discussed in section \ref{sec:prop}. Finally, the numerical examples are presented in section \ref{sec:examp}.

\section{Urysohn operator for control systems}
\label{sec:theory}

The general form of the Urysohn model is a multiple-input multiple-output (MIMO) model \cite{Krylov1979a}. However, this paper focuses mainly on single-input single-output (SISO) simplification of the model. The extension of the proposed modelling and identification techniques to the MIMO case is considered separately in section \ref{sec:multidim}. The relation of the Urysohn model to the linear, the Hammerstein, the Volterra series and the NARMAX models is discussed in appendix \ref{sec:relAll} to indicate the place of the Urysohn model in the hierarchy of control systems' models. 

\subsection{Continuous-time form of the model}

The continuous-time Urysohn operator is an integral operator, which transforms function $x\left(t\right)$ to function $y\left(t\right)$ in the following way \cite{Krylov1979a,Krylov1984a}:
\begin{equation}
  y\left( t \right) = \int_{0}^{T} V\left( s, x\left(t - s\right) \right) \dif s ,
  \label{eq:UryshonCont}
\end{equation}
where $x : \left[ -T , +\infty \right) \to \left[x_\mathrm{min}, x_\mathrm{max}\right]$, $y : \left[ 0 , +\infty \right) \to \mathbb{R}$, $V : \left[ 0 , T \right] \times \left[x_\mathrm{min}, x_\mathrm{max}\right] \to \mathbb{R}$ are continuous almost everywhere functions, $t \geq 0$, $T \geq 0$ and $x_\mathrm{min}, x_\mathrm{max} \in \mathbb{R}$. Function $V\left( s, x \right)$ is the kernel of the continuous-time Urysohn operator\footnote{In the original literature, the operator with an infinite memory is considered, $y\left( t \right) = \int_{-\infty}^{t} V\left( t-\xi, x\left(\xi\right) \right) \dif \xi$. In this paper, following \cite{Poluektov1990}, the operator with a finite memory is considered, i.e. it is assumed that $y\left( t \right)$ is defined completely by $x\left(\xi\right)$, where $t-T \leq \xi \leq t$. By substitution $s = t - \xi$, equation \eqref{eq:UryshonCont} is obtained.}.

In the case of control systems, $x\left(t\right)$ is the time-dependent input of the control system, $y\left(t\right)$ is the time-dependent output of the control system. Argument $\left(t - s\right)$ in equation \eqref{eq:UryshonCont} describes the causality between action $x\left(t\right)$ and reaction $y\left(t\right)$ of the object. Parameter $T$ is the time interval, which is sufficiently large for each $y\left(t\right)$ to be determined by the input within the interval between $\left(t-T\right)$ and $t$. 

\subsection{Discrete-time form of the model}
\label{sec:discrUrysohn}

The discrete-time Urysohn operator is given by \cite{Krylov1979,Krylov1984a}:
\begin{equation}
  y_i = \sum_{j=1}^{m} g_j\left( x_{i-j+1} \right) , \quad i \in \mathbb{N} ,
  \label{eq:UryshonDiscrOrig}
\end{equation}
where $x_i \in \left[x_\mathrm{min}, x_\mathrm{max}\right]$ is the series of input values, $y_i \in \mathbb{R}$ is the series of output values, $m$ is the memory depth of the operator and $g_j: \left[x_\mathrm{min}, x_\mathrm{max}\right] \to \mathbb{R}$ are continuous almost everywhere functions. The set of functions $g_j$ can be called the kernel of the discrete-time Urysohn model. It can be seen that equation \eqref{eq:UryshonDiscrOrig} results from a numerical quadrature\footnote{For example, by defining $x_i = x\left(i \Delta t\right)$ and $y_i = y\left(i \Delta t\right)$, where $\Delta t = T/\left(m-1\right)$ is the time step, and by using the composite trapezoidal rule, equation \eqref{eq:UryshonCont} becomes
\begin{multline*}
  y_i = y\left(i \Delta t\right) = \int_{0}^{T} V\left( s, x\left(i \Delta t - s\right) \right) \dif s \approx
    \frac{1}{2}\Delta t \left( V\left( 0, x\left(i \Delta t\right)\right) +
    2 V\left( \Delta t, x\left(\left(i-1\right) \Delta t\right) \right) + \ldots + \vphantom{V}\right. \\
    \left.\vphantom{V} + 2 V\left( \left(m-2\right) \Delta t, x\left(\left(i-m+2\right) \Delta t\right) \right) +
    V\left( \left(m-1\right)\Delta t, x\left(\left(i-m+1\right) \Delta t\right) \right) \right) .
\end{multline*}
Equation \eqref{eq:UryshonDiscrOrig} is obtained by defining $g_1\left( x \right) = \frac{\Delta t}{2} V\left( 0, x \right), \; \ldots, \; g_m\left( x \right) = \frac{\Delta t}{2} V\left( \left(m-1\right)\Delta t, x \right)$.} of equation \eqref{eq:UryshonCont}.

The discrete-time Urysohn model, equation \eqref{eq:UryshonDiscrOrig}, contains non-linear functions $g_j$ that must be represented in some parametric form before the model can be identified and used to reproduce the input-output relation of a control system. For example, it is possible to represent them by known polynomials with some coefficients \cite{Makarov1994,Makarov2012}. In \cite{Poluektov1990}, it has been proposed to take functions $g_j$ to be piecewise constant --- the input range $\left[ x_\mathrm{min}, x_\mathrm{max} \right]$ is divided into $n$ intervals and functions $g_j$ are constant within each interval. Such model can be conveniently rewritten using the following matrix notation\footnote{Here, the $j$-th row of matrix $U$ contains all $n$ values that piecewise-constant function $g_j$ can take.}:
\begin{align}
  &y_i = \sum_{j=1}^{m} U\left[ j, k_{i-j+1} \right] ,
  \label{eq:UryshonDiscr} \\
  &k_i = 1 + \operatorname{round}\left( \frac{ \left(n-1\right) \left(x_i-x_\mathrm{min}\right) }
    { x_\mathrm{max} - x_\mathrm{min} } \right) ,
  \label{eq:UryshonDiscrContr}
\end{align}
where $U$ is the matrix with indices shown in $\left[ \cdot, \cdot \right]$ and operator $\operatorname{round}\left(\cdot\right)$ is the rounding to the nearest integer. Matrix $U$ is referred to as the Urysohn matrix in the rest of the paper.

Model \eqref{eq:UryshonDiscr}-\eqref{eq:UryshonDiscrContr} can be called the quantised-input discrete-time Urysohn model, as equation \eqref{eq:UryshonDiscrContr} performs the quantisation of the input into $n$ levels, before it is used for the calculation of the output. Here, integer $k_i$ is the quantised input. In the rest of the paper, this model is referred to as the discrete Urysohn model for simplicity, and the operator in equation \eqref{eq:UryshonDiscr} that transforms sequence $k_i$ to sequence $y_i$ is referred to as the discrete Urysohn operator. In section \ref{sec:nonQuantInput}, a more general case of a piecewise-linear representation of $g_j$, i.e. the case of a non-quantised input, is introduced.

It should be noted that model \eqref{eq:UryshonDiscr}-\eqref{eq:UryshonDiscrContr} is not a classical lookup table linking the input and the output of a system, but rather a grid (with matrix $U$ containing the grid point values) with a certain rule for selecting a subset of elements from it, sum of which forms the output. Furthermore, for the case of a non-quantised input, section \ref{sec:nonQuantInput}, certain weights are introduced in the sum. In the case of a classical lookup table, each distinct input sequence is associated with a single distinct output taken directly from a lookup table. Such models are also referred to as non-parametric models \cite{Schoukens2019}, and further discussion of these models can be found in textbooks, e.g. \cite{Fan1996}.

\section{Identification of the discrete Urysohn operator}
\label{sec:ident}

The identification problem for the discrete Urysohn operator consists in finding the unknown elements of matrix $U$ using known input and output sequences. The Urysohn matrix, however, has an important property --- any given input and output sequences do not uniquely determine matrix $U$. More specifically, the Urysohn matrix of size $m \times n$ contains $\left(m - 1\right)$ elements that can be selected arbitrarily. This property can be formulated into the following theorem.

\begin{theorem}
\label{th:ident}
For any given quantised input sequence $k_i$ and output sequence $y_i$ of an Urysohn system \eqref{eq:UryshonDiscr}, when $m > 1$, there are infinitely many Urysohn matrices, for which the input sequence is converted exactly to the output sequence. Moreover, when $\left(m - 1\right)$ elements of the Urysohn matrix are prescribed, selected such that not more than one element from each row of the Urysohn matrix is prescribed, there is a unique set of remaining $\left(m n - m + 1\right)$ elements, such that the discrete Urysohn operator exactly converts the input sequence to the output sequence. In this case, these remaining $\left(m n - m + 1\right)$ elements linearly depend on the values of the prescribed $\left(m - 1\right)$ elements.
\end{theorem}

\begin{proof}
See appendix \ref{sec:rank}.
\end{proof}

\begin{remark}
It must be emphasised that the theorem is formulated for the input and output sequences of an Urysohn system \eqref{eq:UryshonDiscr}, i.e. the input is quantised and the output is formed by the discrete Urysohn operator. Thus, there is at least one solution of the identification problem.
\end{remark}

\begin{remark}
In the theorem, it is implied that the input and the output sequences have the same length; however, first $\left(m - 1\right)$ elements of the output sequence are not defined as the output of the discrete Urysohn operator.
\end{remark}

To reduce the infinity of possible solutions of the identification problem to a unique solution, additional constraints must be introduced. One of such possible constraints is the minimum Frobenius norm of the Urysohn matrix. In section \ref{sec:identIter}, the iterative identification method is proposed, which converges to the unique solution with such minimum norm. 

\subsection{The iterative identification method}
\label{sec:identIter}

The proposed algorithm for identifying the discrete Urysohn operator is based on the Kaczmarz algorithm \cite{Kaczmarz1937,Tewarson1969} for solving linear systems of equations. The key step that allows using the Kaczmarz algorithm is the assembly of the linear system of equations with respect to parameters of the discrete Urysohn operator. This is done by performing a logical operation on the input of non-linear system \eqref{eq:UryshonDiscr}-\eqref{eq:UryshonDiscrContr}. The details are given in section \ref{sec:identKacz}.

The proposed algorithm can be summarised as follows:
\begin{enumerate}
  \item Assume an initial approximation of matrix $U$. It can be arbitrary, including the all-zero matrix.
  \item Start with $i = m$.
  \item Calculate model output $\hat{y}_i$ based on actual input $x_i$ and the current approximation of matrix $U$ according to equations \eqref{eq:UryshonDiscr}-\eqref{eq:UryshonDiscrContr}. \label{item:iterS}
  \item Calculate difference $D = y_i - \hat{y}_i$, where $y_i$ is the actual recorded output and $\hat{y}_i$ is the model output.
  \item Modify matrix $U$, such that $\alpha D / m$ is added to each element that was involved in the calculation of $\hat{y}_i$. \label{item:iterE}
  \item Increase index $i$ by $1$ and repeat steps \ref{item:iterS}-\ref{item:iterE} until $D$ becomes sufficiently small for sufficiently large number of iterations consecutively.
\end{enumerate}
In step \ref{item:iterE} of the algorithm, $\alpha D / m$ is added to elements $\left[ j, k_{i-j+1} \right]$ of matrix $U$, where $j = 1,\ldots,m$. Parameter $\alpha$ is the stabilisation parameter from interval $\left(0, 1\right]$ for suppressing the noise. For near exact data, $\alpha$ can be $1$, while for a very noisy data, it must be relatively low. Parameter $\alpha$ is further discussed in sections \ref{sec:resultsNoisyO} and \ref{sec:resultsNoisyIO}, where numerical examples are provided.

The theorem on the uniqueness of the solution is formulated below. Furthermore, an independent proof of the theorem (without using the relation to the Kaczmarz algorithm) is given.

\begin{theorem}
\label{th:conv}
For exact input-output data of an Urysohn system, when the input sequence covers all possible inputs within the desired range and the values of the input sequence from $\left(i-m+1\right)$ to $i$ almost always change\footnote{When fragment of the input sequence from $\left(i-m+1\right)$ to $i$ does not change with the iteration number, matrix $U$ is not updated, i.e. the repeated sequences are just ignored by the algorithm.} with iteration number $i$, the proposed algorithm presented above converges to a unique solution. If the initial approximation for matrix $U$ is the all-zero matrix, the proposed algorithm converges to matrix $U$ with the minimum Frobenius norm.
\end{theorem}

\begin{proof}
See appendix \ref{sec:identConvProof}.
\end{proof}

The entire block of adjustment operations is computationally inexpensive and can be applied as a real-time process in an interval between automatic reading of the measurements of the input and the output of the physical control system. In an automatic identification, it is easy to trace the number of times each element of matrix $U$ has been modified, which can be an important information for assessing whether the algorithm has converged.

\subsection{Relation to the Kaczmarz iterative method}
\label{sec:identKacz}

The Kaczmarz iterative method for solving linear system of equations $A X = B$, where $A$ is a known matrix, $B$ is a known vector-column and $X$ is an unknown vector-column, is given by the following sequence \cite{Kaczmarz1937,Tewarson1969}:
\begin{equation}
  X^{i+1} = X^i + \frac{B_p - A_p X^i}{\left| A_p \right|^2} {A_p}^\mathrm{T} ,
  \label{eq:KaczGen}
\end{equation}
where $A_p$ is the $p$-th row of $A$, $B_p$ is the $p$-th element of $B$ and $X^i$ is the approximation of the solution at step $i$. Integer variable $p$ changes with the iteration number.

In the case of the discrete Uryshon operator, each value of the output is the sum of the specifically selected subset of elements of matrix $U$. This allows building a system of linear algebraic equations for the identification of matrix $U$ by rearranging its elements into an unknown vector-column
\begin{equation*}
  Z = \begin{bmatrix}
    U_{11} & \ldots & U_{1n} & U_{21} & \ldots & U_{2n} & \ldots & U_{m1} & \ldots & U_{mn}
  \end{bmatrix}\transp .
\end{equation*}
Recorded output sequence $y_i$ is also rearranged into a vector-column
\begin{equation*}
  \tilde{Y} = \begin{bmatrix}
    y_{m} & y_{m+1} & \ldots & y_{m+N}
  \end{bmatrix}\transp ,
\end{equation*}
where $N$ is the number of elements in the recorded output sequence. Matrix $\tilde{M}$ is introduced, elements of which are determined using known $k_i$ as
\begin{equation*}
  \tilde{M}_{iq} = \begin{cases}
    1, & \text{if } q = n\left(j-1\right) + k_{i-j+1} \text{ where } j\in\left\lbrace 1, 2, \ldots, m \right\rbrace \\
    0, & \text{otherwise.}
  \end{cases}
\end{equation*}
In this case, the discrete Urysohn system \eqref{eq:UryshonDiscr} can be represented by the following system of algebraic equations:
\begin{equation}
  \tilde{M} Z = \tilde{Y} .
  \label{eq:UrysohnSLAE}
\end{equation}
Thus, in order to find unknown $Z$, the system of linear equations \eqref{eq:UrysohnSLAE} must be solved. The rank of matrix $\tilde{M}$ is $\left(m n - m + 1\right)$ or less, due to theorem \ref{th:ident}. The rank of $\tilde{M}$ is strictly less than $\left(m n - m + 1\right)$ when the input sequence does not cover all possible input values.

It can be seen that by applying the Kaczmarz method to system \eqref{eq:UrysohnSLAE} and by introducing a multiplier $\alpha$ in equation \eqref{eq:KaczGen}, the algorithm of section \ref{sec:identIter} is obtained. The norm of the each row of $\tilde{M}$ is $\sqrt{m}$ and, at each iterative step, only those elements of $Z$ are modified that were involved in the calculation of the corresponding element of $\tilde{Y}$.

The Kaczmarz method is also sometimes called the projection descent method \cite{Faddeev1981}, which results from its geometrical interpretation. The solution of a linear system of equations can be interpreted as finding an intersection point of hyperplanes in a multidimensional space. Initially, an arbitrary point is taken and, at each iteration, the point is projected onto a different hyperplane. Each projection operation brings the point closer to the solution.

The actual convergence rate depends significantly on the angles between the hyperplanes \cite{Censor2009}. For close to orthogonal set of hyperplanes, the convergence is relatively fast, while for hyperplanes intersecting at sharp angles, the convergence is relatively slow. The convergence of the Kaczmarz method in application to system \eqref{eq:UrysohnSLAE} is relatively fast due to rows of matrix $\tilde{M}$ being either orthogonal or relatively close to being orthogonal, since, for each new input/output element, most non-zero input values of matrix $\tilde{M}$ are expected to be shifted to neighbouring positions in the matrix\footnote{This can be illustrated by the following example. Assume $m = 3$, $n = 3$ and $k_1 = 1$,  $k_2 = 2$,  $k_3 = 1$,  $k_4 = 3$. Then rows of $\tilde{M}$ corresponding to $y_3$ and $y_4$ are
\begin{align*}
  &\begin{bmatrix}
    1 & 0 & 0 & 0 & 1 & 0 & 1 & 0 & 0
  \end{bmatrix} , \\
  &\begin{bmatrix}
    0 & 0 & 1 & 1 & 0 & 0 & 0 & 1 & 0
  \end{bmatrix} ,
\end{align*}
respectively, which are orthogonal.
}.

The identification of a discrete Urysohn system has already been considered in \cite{Poluektov1990}, where the identification of the operator has been performed by a direct solution of system \eqref{eq:UrysohnSLAE} using the Tikhonov regularisation. To avoid dealing with the degenerate matrix, $\left(m - 1\right)$ values of the Urysohn matrix have been fixed. The disadvantage of such method is that in the case of noisy data, elements of matrix $\tilde{M}$ are slightly misplaced and may be located at adjacent positions, which complicates the solution of system \eqref{eq:UrysohnSLAE} even when the regularisation is utilised. Furthermore, as the method of \cite{Poluektov1990} operates with the fully assembled matrices, it requires significantly larger memory and cannot be used in real time for model identification, in contrast to the method proposed above.

It should be noted that the proposed algorithm is similar to the well-known Least Mean Squares (LMS) algorithm \cite{Widrow1960,Haykin2014} for identification of linear control systems. The similarity comes from the fact that the Normalised Least Mean Squares (NLMS) algorithm is the Kaczmarz algorithm\footnote{With a possible difference that NLMS can include a variable step size \cite{Haykin2014}.} applied to matrix form of the linear control system \cite{Anjum2015}. The major difference between the proposed algorithm and the LMS/NLMS algorithms is a set of logical operations preceding the model update. As seen form the structure of the discrete Urysohn model, the input sequence defines the addresses of a subset of elements of the identified grid structure. The elements of the grid structure, in turn, form the output of the model. This contrasts with the LMS/NLMS algorithms for linear systems, where a linear combination of the elements of the input is equal to the output, and the coefficients of this linear combination are identified.

\section{Some generalisations and properties of the discrete Urysohn operator}
\label{sec:prop}

\subsection{Generalisation for non-quantised input}
\label{sec:nonQuantInput}

The discrete Urysohn model, which has been introduced in section \ref{sec:discrUrysohn}, involves the quantisation of the input. Model \eqref{eq:UryshonDiscr}-\eqref{eq:UryshonDiscrContr} is based on piecewise-constant representation of $g_j$ in equation \eqref{eq:UryshonDiscrOrig}. A piecewise-linear representation of $g_j$, on the other hand, leads to a more general model, where the quantisation of the input is not required. For simplicity of the presentation, the input range $\left[ x_\mathrm{min}, x_\mathrm{max} \right]$ is divided into $n$ equal intervals and functions $g_j$ are taken to be linear within each interval. In this case, the discrete-time Urysohn operator is still represented by a matrix\footnote{In the case of such piecewise-linear representation of $g_j$, these functions are fully defined by the nodal values of $g_j$, i.e. the values of the functions where the slope changes. In the proposed generalisation, the $j$-th row of matrix $U$ contains all nodal values of $g_j$.}; however, the rule for calculating the output changes, as is shown in this section.

First, a rescaling of the input is introduced --- an additional variable $b_i$ is constructed in the following way:
\begin{equation}
  b_i = 1 + \left( n-1 \right) \frac{x_i - x_\mathrm{min}}{x_\mathrm{max} - x_\mathrm{min}} .
\end{equation}
This variable and takes all real values from interval $\left[ 1, n \right]$, as $x_i$ takes all real values from interval $\left[ x_\mathrm{min}, x_\mathrm{max} \right]$. Next, rounding to the nearest integer values is introduced:
\begin{equation}
  k_i^\mathrm{L} = \left\lfloor b_i \right\rfloor , \quad
  k_i^\mathrm{R} = \left\lceil b_i \right\rceil ,
\end{equation}
where $\lfloor \cdot \rfloor$ and $\lceil \cdot \rceil$ are the floor and the ceiling functions, respectively. These integers are needed to address the elements of the Urysohn matrix. Finally, the generalised form of the operator can be introduced:
\begin{equation}
  y_i = \sum_{j=1}^{m} \left( \left( 1 - \psi_{i-j+1} \right) U\left[ j, k^\mathrm{L}_{i-j+1} \right] +
  \psi_{i-j+1} U\left[ j, k^\mathrm{R}_{i-j+1} \right] \right) , \quad
  \psi_i = b_i - k_i^\mathrm{L} .
\end{equation}
It can be seen that in the above representation, each term of the sum changes piecewise-linearly as a function of $b_i$, hence as a function of $x_i$, and the nodal values (i.e. the values at the points where the slope changes) are taken from the $j$-th row of matrix $U$.

The identification procedure for the non-quantised case also changes. Following the general formula of the Kaczmarz method, equation \eqref{eq:KaczGen}, the following norm is introduced:
\begin{equation}
  \chi_i = \sum_{j=1}^{m} \left( \left( 1 - \psi_{i-j+1} \right)^2 + {\psi_{i-j+1}}^2 \right) .
\end{equation}
In the non-quantised case, only step \ref{item:iterE} of the algorithm of section \ref{sec:identIter} changes. Now, at this step, $\alpha D \left( 1 - \psi_{i-j+1} \right) / \chi_i$ is added to elements $\left[ j, k^\mathrm{L}_{i-j+1} \right]$ of matrix $U$ and  $\alpha D \psi_{i-j+1} / \chi_i$ is added to elements $\left[ j, k^\mathrm{R}_{i-j+1} \right]$ of matrix $U$. It should be emphasised that each element of matrix $U$ is modified by a different value, which contains the corresponding weight.

\subsection{Multiple inputs and multiple outputs}
\label{sec:multidim}

The general form of the Urysohn model is a MIMO model \cite{Krylov1979a}. The approach of section \ref{sec:nonQuantInput} can easily be generalised to the MIMO case. For simplicity of the presentation, first, a two-input single-output system is considered. The continuous-time form of the model is given by \cite{Krylov1979a}
\begin{equation}
  z\left( t \right) = \int_{0}^{T} V\left( s, x\left(t - s\right), y\left(t - s\right) \right) \dif s ,
  \label{eq:UryshonMltCont}
\end{equation}
where $x\left(t\right)$ and $y\left(t\right)$ are the inputs and $z\left(t\right)$ is the output. The corresponding discrete-time form is given by \cite{Krylov1979}
\begin{equation}
  z_i = \sum_{j=1}^{m} g_j\left( x_{i-j+1}, y_{i-j+1} \right) ,
\end{equation}
where $x_i$ and $y_i$ are the input sequences and $z_i$ is the output sequence. Following the idea of piecewise-linear representation of functions for the single-input case, as $g_j$ are now functions of two variables, a piecewise-bilinear representation can be used. Following similar steps as in section \ref{sec:nonQuantInput}, results in the following model:
\begin{align}
  &z_i = \sum_{j=1}^{m} \left( \left( 1 - \psi_{i-j+1} \right) \left( 1 - \phi_{i-j+1} \right)
    U\left[ j, k^\mathrm{L}_{i-j+1}, k^\mathrm{D}_{i-j+1} \right] +
    \vphantom{U}\right.\nonumber \\ &\qquad\qquad\left.\vphantom{U} 
    \psi_{i-j+1} \left( 1 - \phi_{i-j+1} \right) U\left[ j, k^\mathrm{R}_{i-j+1}, k^\mathrm{D}_{i-j+1} \right] + 
    \vphantom{U}\right.\nonumber \\ &\qquad\qquad\left.\vphantom{U} 
    \left( 1 - \psi_{i-j+1} \right) \phi_{i-j+1} U\left[ j, k^\mathrm{L}_{i-j+1}, k^\mathrm{U}_{i-j+1} \right] +
    \vphantom{U}\right.\nonumber \\ &\qquad\qquad\left.\vphantom{U} 
    \psi_{i-j+1} \phi_{i-j+1} U\left[ j, k^\mathrm{R}_{i-j+1}, k^\mathrm{U}_{i-j+1} \right] \right) , \\
  &\psi_i = b_i - k_i^\mathrm{L} , \quad 
    \phi_i = c_i - k_i^\mathrm{D} , \nonumber \\
  &k_i^\mathrm{L} = \left\lfloor b_i \right\rfloor , \quad
    k_i^\mathrm{R} = \left\lceil b_i \right\rceil , \quad
    k_i^\mathrm{D} = \left\lfloor c_i \right\rfloor , \quad
    k_i^\mathrm{U} = \left\lceil c_i \right\rceil , \\
  &b_i = 1 + \left( n-1 \right) \frac{x_i - x_\mathrm{min}}{x_\mathrm{max} - x_\mathrm{min}} , \quad
    c_i = 1 + \left( n-1 \right) \frac{y_i - y_\mathrm{min}}{y_\mathrm{max} - y_\mathrm{min}} , 
\end{align}
where $U$ is now the three-dimensional matrix with indices shown in $\left[ \cdot, \cdot, \cdot \right]$; variables $b_i$ and $c_i$ are scaled inputs $x_i$ and $y_i$; integers $k_i^\mathrm{L}$ and $k_i^\mathrm{R}$ are rounding of variable $b_i$ down and up, respectively; integers $k_i^\mathrm{D}$ and $k_i^\mathrm{U}$ are rounding of variable $c_i$ down and up, respectively. These integers are used to address elements of matrix $U$.

To perform the identification procedure, following equation \eqref{eq:KaczGen} and steps of section \ref{sec:identKacz}, the following norm is introduced:
\begin{align*}
  &\chi_i = \sum_{j=1}^{m} \left( \left( 1 - \psi_{i-j+1} \right)^2 \left( 1 - \phi_{i-j+1} \right)^2 +
    {\psi_{i-j+1}}^2 \left( 1 - \phi_{i-j+1} \right)^2 + 
    \vphantom{U}\right. \\ &\qquad\qquad\left.\vphantom{U}
    \left( 1 - \psi_{i-j+1} \right)^2 {\phi_{i-j+1}}^2 +
    {\psi_{i-j+1}}^2 {\phi_{i-j+1}}^2 \right) .
\end{align*}
Again, only step \ref{item:iterE} of the algorithm of section \ref{sec:identIter} changes. Now, at this step, $\alpha D \left( 1 - \psi_{i-j+1} \right) \left( 1 - \phi_{i-j+1} \right) / \chi_i$ is added to elements $\left[ j, k^\mathrm{L}_{i-j+1}, k^\mathrm{D}_{i-j+1} \right]$ of matrix $U$, $\alpha D \psi_{i-j+1} \left( 1 - \phi_{i-j+1} \right) / \chi_i$ is added to elements $\left[ j, k^\mathrm{R}_{i-j+1}, k^\mathrm{D}_{i-j+1} \right]$ of matrix $U$, etc.

It easy to see that the presented above approach can be trivially generalised for multiple-input single-output (MISO) system by representing corresponding functions $g_i$ using piecewise-multilinear (e.g. piecewise-trilinear for three inputs) representation. The MIMO case is obtained from the MISO case by constructing an individual model for each output. The outputs may be correlated, but they are independent by the definition of the Urysohn model, and they depend only on the inputs.  

\subsection{Describability of a system by the discrete Urysohn operator}
\label{sec:descrCrit}

The integral models have certain advantages over differential ones. For example, the convolution-type linear integral equation can describe an object of any order with a pure delay. Thus, it is useful to provide a criteria for a system to be describable by the Urysohn operator.

\begin{theorem}
\label{th:syst}
For a system with a quantised input to be describable by the discrete Urysohn operator, it is necessary and sufficient that (a) the system has a finite memory (b) the system has additivity for input sequences that are individual impulses, which do not coincide in time.
\end{theorem}

\begin{proof}
See appendix \ref{sec:critProof}.
\end{proof}

\begin{remark}
Since the case of a quantised input is considered, an individual impulse means that all elements of the input sequence are equal to $x_\mathrm{min}$ except one element, which can take an arbitrary quantised value.
\end{remark}

\begin{remark}
The condition that the impulses do not coincide in time is essential. Condition (b) is not a linearity condition. If a system has additivity for impulses that coincide in time, this means that the system is linear with respect to inputs. In the general form, the Urysohn systems are non-linear and, therefore, do not have additivity for arbitrary inputs.
\end{remark}

In the case of a non-quantised input, following equation \eqref{eq:UryshonDiscrOrig}, the Urysohn systems are such systems, where the output is a sum of non-linear functions of elements of the input sequence. Systems, where the output depends on products of time-shifted elements of the input sequence, cannot be exactly represented by the Urysohn model.

In the case of real applications, it is easy to judge the applicability of the discrete Urysohn model to a system. In the case when the above criteria are not strictly fulfilled, the error (the deviation from the criteria) can give an indication of how accurate the discrete Urysohn model can describe the system. However, such investigation requires an ability to impose arbitrary input signals (impulses), which is not always possible. The identification algorithm of section \ref{sec:identIter}, on the other hand, does not require this and works with the observed signals. This means that from a practical point of view, in some cases, it can be easier to apply the identification procedure and judge the applicability of the Urysohn model, based on the accuracy of the representation rather that first check whether the system is of the Urysohn type. In this case, it would also be desirable to show that a simpler model, such as the Hammerstein model, cannot represent the system, i.e. the error is larger if the simpler model is used. Otherwise, the simpler model is always preferable.

\subsection{Partially identified Urysohn operator}
\label{sec:partOper}

The discrete Urysohn operator has a unique property --- even a partially identified operator is still useful and can be utilised. The operator will be identified partially when input values, which are used for the identification, do not cover the entire range from $x_\mathrm{min}$ to $x_\mathrm{max}$. Thus the notion of ``identification range'' can be introduced as the range, within which the input varies during the identification. In this case, as evident from the identification procedure, elements of the Urysohn matrix corresponding to the input outside of this range will not be updated and will remain to be initial guesses. However, the ``middle'' elements of the Urysohn matrix, which correspond to input values within the identification range, will be identified. 

Using the proposed identification procedure, it is easy to introduce the counter for every element of $U$, which stores the number of times the element has been updated. Such counter can be useful to determine the identification range and to estimate the accuracy of the identified elements, since the error decreases with the number of updates (due to the convergence of the identification procedure).

When such partially identified operator is applied to a different input sequence, it will still produce a reliable output when the input is within the identification range. When the input temporarily takes values outside of the identification range, the model does not produce an output; however, when the input comes back into the identification range and stays there for a time period, which is greater than the system memory, the model again starts producing a valid output.

\section{Numerical examples}
\label{sec:examp}

The goal of this section is to demonstrate the descriptive capabilities of the discrete Urysohn operator. Two series of test studies are performed. For the first set of tests, a non-linear controllable mechanical object, the dynamic behaviour of which is described by a non-linear differential equation, is considered. Exact input and output sequences are generated using the numerical solution of the differential equation. This allows systematic studying of the performance of the identification algorithm. The second set of tests is conducted using real experimental data and aims  at showing the readers that the discrete Urysohn operator performs well in real-world scenarios.

\subsection{Studied system}
\label{sec:system}

\begin{figure}
  \begin{center}
    \includegraphics{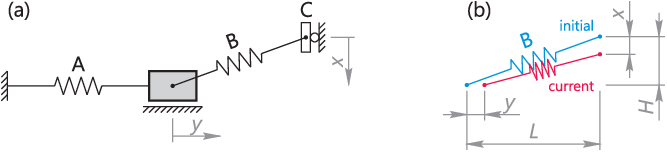}
  \end{center}
  \caption{A schematic representation of the considered mechanical system (a) and a schematic representation of the change of the geometry of spring B due to the movement of the object (b).}
  \label{fig:obj}
\end{figure}

The considered mechanical system is shown in figure \ref{fig:obj}a. A bulky object is allowed to move in the horizontal direction; the movement is affected by a friction force. The object is connected to a clamp by a horizontal spring. A second spring is connected by a hinge to the object and by another hinge to a platform, which can move vertically. The vertical displacement of the platform is the input of the system (the control), while the horizontal displacement of the object is the output of the system (the observable state variable).

Although the components of the system (the springs) are linear, the static input-output response of the system (the steady-state characteristic) is significantly non-linear due to the geometry of the system. In the static case, $y$ decreases with the increase of $x$ up to a point when platform C is parallel to the object. Afterwards, $y$ increases with the further increase of $x$.

The mechanical system is described by the following differential equation:
\begin{align}
  &\frac{\dif^2 y}{\dif t^2} = -2\zeta\omega\frac{\dif y}{\dif t} + f\left(y,x\right) ,
  \label{eq:state} \\
  &f\left(y,x\right) = -\omega^2 y - \omega^2 \left( \sqrt{L^2+H^2} -
    \sqrt{\left(L-y\right)^2+\left(H-x\right)^2} \right) \frac{L-y}{\sqrt{\left(L-y\right)^2+\left(H-x\right)^2}} ,
    \\
  &\left. y \right|_{t=0} = y_0 , \quad
  \left. \frac{\dif y}{\dif t} \right|_{t=0} = v_0 , 
\end{align}
where $y$ and $x$ are the state variable of the system and the control function, respectively; $2\zeta\omega$ is the friction coefficient divided by the mass of the object; $\omega^2$ is the stiffness of the springs divided by the mass of the object; $L$ and $H$ are the initial horizontal and vertical distances, respectively, between the centres of the object and the platform. The expression for the total force scaled by the mass, $f\left(y,x\right)$, results from the projection of the force in spring B onto the horizontal direction. The change of the geometry of spring B is shown in figure \ref{fig:obj}b.

The numerical solution of equation \eqref{eq:state} can be obtained using the Verlet method. It is easy to verify that the discretisation of equation \eqref{eq:state} using the Verlet method results in
\begin{align}
  &y_{i+1} = \left( 2 y_i - y_{i-1} \left( 1 - 2\zeta\omega\frac{\Delta t}{2} \right) +
    f\left(y_i,x_i\right) \Delta t ^2 \right) \left( 1 + 2\zeta\omega\frac{\Delta t}{2} \right)^{-1},
  \label{eq:stateDiscr} \\
  &y_1 = y_0 + \Delta t v_0 \left( 1 - 2\zeta\omega\frac{\Delta t}{2} \right) +
    f\left(y_0,x_0\right) \frac{\Delta t ^2}{2} ,
  \label{eq:stateDiscrStart}
\end{align}
where $\Delta t$ is the time step and subscript $i$ indicates that a quantity is taken a time step $i$.

For the purpose of this paper, all quantities in the equations are taken to be dimensionless. Since it is possible to perform spatial and temporal scaling of the system, which does not affect the qualitative behaviour of the system, parameters $\omega$ and $L$ can be chosen arbitrarily. Parameter $H/L$ controls the degree of non-linearity, while parameter $\zeta$ controls the oscillatory nature of the system. Numerical parameter $\Delta t$ must be chosen in such way that the numerical results maintain sufficient accuracy and, for oscillatory systems, it is usually selected as a fraction of the period of undamped oscillations. The following values of the parameters are taken: $\omega = 1$, $\zeta = 1$, $L = 1$, $H = 0.5$, $\Delta t = 2 \pi / 128$. Initially, the system is at rest: $y_0 = 0$, $v_0 = 0$.

\subsection{Results of identification}
\label{sec:results}

\subsubsection{Discrete control function}
\label{sec:resultsDiscr}

Since the discrete Urysohn operator requires certain discretisation of the control function and also certain discretisation in time, the simplest case for the identification is when the control function takes only discrete values and is constant for periods of $\Delta \tau$. Thus, the following control function is considered:
\begin{align}
  &x\left(t\right) = \left(k-1\right) \Delta x \quad\mbox{for}\quad \left( j-1 \right)
    \Delta \tau < t \leq j \Delta \tau ,
  \label{eq:discrContr} \\
  &k \in \left\lbrace 1,2,\ldots,n \right\rbrace, \quad j \in \left\lbrace 1,2,\ldots,Q \right\rbrace,
    \quad Q = \operatorname{round}\left(t_\mathrm{max} / \Delta \tau\right) . \nonumber
\end{align}
Values $\Delta \tau = 2 \pi / 8$ and $\Delta x = 0.1$ are selected for the numerical experiments based on the dynamic properties of the system. The number of rows and columns of the Urysohn matrix is selected to be $m=8$ and $n = 11$, respectively, which gives the maximum value of the control function $x_\mathrm{max} = 1$.

Any proper benchmarking of any identification procedure requires two completely independent input-output datasets --- the first dataset for the identification of the model, the second dataset for the validation and error calculation. Thus, the model is validated on the unseen data. This strategy is strictly followed in all examples of this paper.

\begin{figure}
  \begin{center}
  \includegraphics{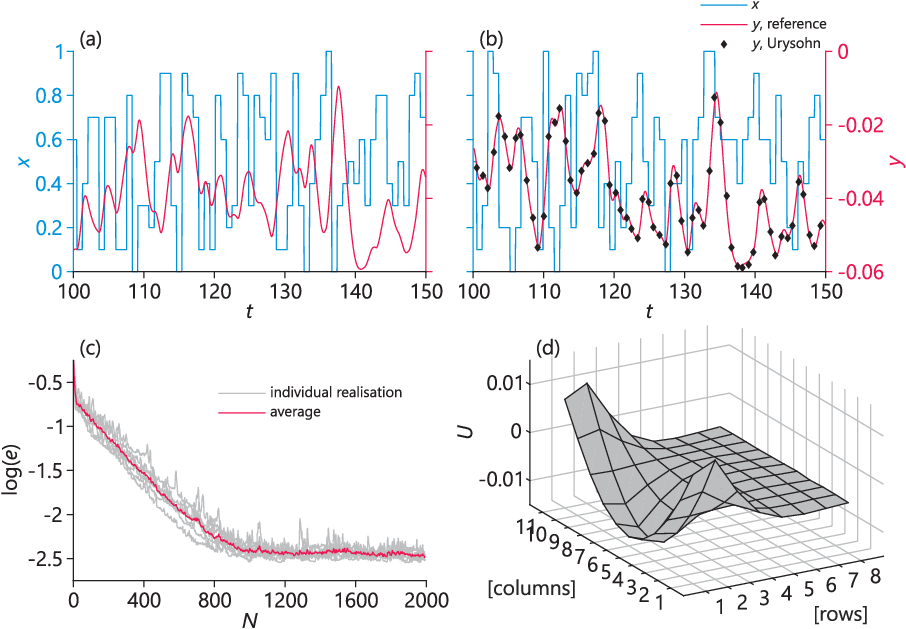}
  \caption{Fragments of the dependence of the input signal and the output of the object on time that were used for the operator identification (a) and validation (b). The Urysohn output is plotted using diamond symbols (b). The dependence of the error of the Urysohn solution on the number of iterations during the identification procedure (c). A visual representation of the Urysohn matrix (d).}
  \label{fig:sol_err_K}
  \end{center}
\end{figure}

The input and the output sequences are generated as follows. Different realisations of the random input signal $x\left(t\right)$ are generated according to equation \eqref{eq:discrContr}. The corresponding outputs of the system $y\left(t\right)$ are calculated using equations \eqref{eq:stateDiscr} and \eqref{eq:stateDiscrStart}. For the purpose of this example, these outputs are considered to be the exact object outputs and are referred to as the reference outputs. The total simulation time is taken to be $t_\mathrm{max} = 10^4$; however, as shown below, much smaller signal length is required for the identification of the Urysohn matrix. A small fragment of the input and the output signals for one of the realisations is shown in figure \ref{fig:sol_err_K}a.

The Urysohn matrices are identified for each realisation of the identification dataset using the algorithm of section \ref{sec:identIter} with $\alpha = 1$. The initial estimates for the Urysohn matrices are the all-zero matrices. Since the computed matrix changes each iteration, it is recorded after each iteration for the subsequent error analysis. An example of the Urysohn matrix, which is obtained after the identification procedure, is illustrated in figure \ref{fig:sol_err_K}d. To validate the obtained Urysohn models and to calculate the errors, the obtained Urysohn matrices are applied to the input signals of the validation dataset and the outputs, which are referred to as the Urysohn outputs, are calculated.

The comparison of the reference output and the Urysohn output is shown in figure \ref{fig:sol_err_K}b. It can be seen that the Urysohn output almost perfectly fits the reference output. The error can be characterised by the scaled $L^1$-norm of the difference between the solutions:
\begin{equation}
  e = \frac{1}{Q y_\mathrm{smax}} \sum_{j=m}^Q \left| \tilde{y}_j - \hat{y}_j \right| , \quad 
  y_\mathrm{smax} = \frac{1}{2}\left( \sqrt{L^2 + H^2} - L \right) ,
  \label{eq:errorL1}
\end{equation}
where $\tilde{y}_j$ is the reference output taken at points $t = j \Delta \tau$ and $\hat{y}_j$ is the Urysohn output. The maximum absolute static displacement $y_\mathrm{smax}$ is introduced to obtain the relative error. The major result of this section is that the average output error of the considered system modelled by the discrete Urysohn operator is $e \approx 0.4\%$ across $8$ different realisations.

Small fluctuations in the output error are related to the randomness of the input sequence. For $\alpha = 1$, which is used in the identification process, the Urysohn matrix resulting from the identification process is highly affected by the last few iterations. To decrease the fluctuations of the output error, parameter $\alpha$ should be decreased.

It is also possible to track the evolution of the error depending on the number of iterations used for the identification of the Urysohn matrix. The Urysohn output and the error are calculated for different Urysohn matrices, where number of iterations $N$ is varied. The error as a function of $N$ is obtained for $8$ different realisations. The results are plotted in figure \ref{fig:sol_err_K}c. The logarithm of the error decreases linearly depending on $N$, until the error reaches a plateau. Since the input of the system is random, the dependence of the error on the number of iterations varies for different realisations.

\subsubsection{Continuous control function}
\label{sec:resultsCont}

The case that has been considered in section \ref{sec:resultsDiscr} is the most simple case in terms of the discretisation of the Urysohn operator. The control function takes only specific discrete values and changes its value at specific moments of time, which are divisible by $\Delta \tau$. In reality, such systems are relatively rare and most systems have a continuous-time control function. In this case, the Urysohn operator must be discretised accordingly, such that a certain accuracy can be achieved.

For a given system, the discretisation of the Urysohn operator is characterised by two parameters --- the discretisation in time, $\Delta \tau$, which is responsible for the number of rows of the Urysohn matrix, and the discretisation of the control function, $\Delta x$, which is responsible for the number of columns of the Urysohn matrix. The desired property of any discrete model is the convergence with respect to the discretisation parameters. Therefore, the aim of this section is to study the accuracy of the discrete Urysohn model depending on $\Delta \tau$ and $\Delta x$.

The convergence of the Urysohn model with respect to the discretisation parameters is verified using the same object as before. At first, the fine-sampled inputs are generated and the corresponding fine-sampled outputs of the system are calculated. For the purpose of this example, these inputs and outputs are regarded as the exact behaviour of the system. After this, the coarse-sampled input and output signals are calculated by local averaging of the fine-sampled signals. The Urysohn matrices are identified based on these coarse-sampled signals. After this, using separate input-output replications, the Urysohn matrices are used to reproduce the output signals. The final step is the comparison of the fine-sampled exact outputs and the Urysohn outputs of the system. These steps are described in detail below.

The control function that corresponds to a random input is considered. Bounded continuous-time function $x\left( t \right)$ is obtained from the following SDE:
\begin{equation}
  \dif x\left( t \right) = G \dif W\left( t \right) ,
  \quad 0 \leq x\left( t \right) \leq 1 ,
  \quad x\left( 0 \right) = 0 ,
  \label{eq:uSDE}
\end{equation}
where $W\left( t \right)$ is the Wiener process and $G$ is the parameter controlling the rate of change of the control function\footnote{In the case when the control function is unbounded and is described by $\dif x\left( t \right) = G \dif W\left( t \right)$, the expected value of the change of $x$ within period $T_0$ can be calculated using
\begin{equation*}
  \Delta X = \int_{-\infty}^{\infty} \left| z \right| f \left( z \right) \dif z ,
\end{equation*}
where $f \left( z \right)$ is the normal distribution function with variance $G^2 T_0$. It is easy to verify that for $T_0 = 2\pi$, $\Delta X = 2G$.}. When $x$ reaches the boundaries, it undergoes the perfect reflection. Numerically, equation \eqref{eq:uSDE} results in the following discrete values of the control function:
\begin{align}
  &p_i = G \sqrt{\Delta t} w_i , \quad i \geq 1 , \\
  &q_i = \sum_{j=1}^i p_j , \\
  &x_i =
  \begin{cases}
    q_i - \left\lfloor q_i \right\rfloor , & \text{if } \left\lfloor q_i \right\rfloor \text{ is even,} \\
    1 - q_i + \left\lfloor q_i \right\rfloor , & \text{if } \left\lfloor q_i \right\rfloor \text{ is odd,} 
  \end{cases}
  \label{eq:uSDE_discr}
\end{align}
where $w_i \sim \mathcal{N}\left( 0, 1 \right)$ are normally distributed random numbers with zero mean and unit variance, $\lfloor \cdot \rfloor$ is the floor function. The structure of equation \eqref{eq:uSDE_discr} takes into consideration reflections from boundaries $0$ and $1$. Parameter $G$ must be selected such that the system reveals its dynamic properties; value $G = 0.05$ is selected for the simulations.

Different realisations of input signal $x_i$ are generated according to equation \eqref{eq:uSDE_discr}. Afterwards, outputs $y_i$ are calculated using equations \eqref{eq:stateDiscr} and \eqref{eq:stateDiscrStart}. The total simulation time is taken to be $t_\mathrm{max} = 10^4$.

To obtain the datasets for the identification and for the validation of the Urysohn matrices, the coarse versions of the input and output signals are calculated for various $\Delta \tau$ and $\Delta x$. Local averaging is used to obtain the coarse-sampled signals:
\begin{align*}
  &x_i^\mathrm{C} = \Delta x \operatorname{round}\left( \frac{1}{\Delta x N_\mathrm{s}}
    \sum_{j = i N_\mathrm{s} - N_\mathrm{s} + 1}^{i N_\mathrm{s}} x_j \right) ,
    \quad\quad N_\mathrm{s} = \frac{ \Delta \tau }{ \Delta t } , \\
  &y_i^\mathrm{C} = \frac{1}{N_\mathrm{s}} \sum_{j = i N_\mathrm{s} - N_\mathrm{s} + 1}^{i N_\mathrm{s}} y_j , 
\end{align*}
where $x_i^\mathrm{C}$ and $y_i^\mathrm{C}$ are the coarse-sampled input and output signals, respectively. The Urysohn matrices are identified based on $x_i^\mathrm{C}$ and $y_i^\mathrm{C}$ for various $\Delta \tau$ and $\Delta x$ using the iterative identification procedure with $\alpha = 1$. The initial estimates for the Urysohn matrices are the all-zero matrices.

The scaled $L^1$-norm of the difference between the exact and the Urysohn outputs, equation \eqref{eq:errorL1}, is again used as the measure for the validation. The results for different discretisation parameters are presented in table \ref{tab:errSto}, where $m = T_0/\Delta \tau$ and $n = 1/\Delta x + 1$ are the number of rows and the number of columns of the Urysohn matrix, respectively, and $T_0 = 2\pi$. It can be seen that the error decreases with the decrease of $\Delta \tau$ and $\Delta x$. The average error that is less than $1\%$ can be achieved for small values of the discretisation parameters.

\begin{figure}
  \begin{center}
    \includegraphics{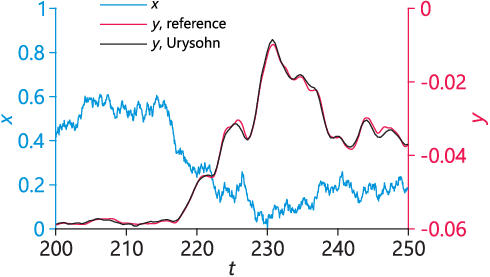}
  \end{center}
  \caption{Fragments of the fine-sampled input and output signals, which were used for the validation. The coarse-sampled Urysohn output for the case of $m=16$ and $n=41$ is plotted using the black line.}
  \label{fig:sol_Ur_sto}
\end{figure}

In figure \ref{fig:sol_Ur_sto}, fragments of the reference inputs and outputs are presented as well as the output of the Urysohn model for the case of $m=16$ and $n=41$. It can be seen that the Urysohn model captures accurately the dynamic behaviour of the system.

\begin{table}
\begin{center}
\begin{tabular}{ | r | c c c c | } 
  \hline
         & $n=11$          & $n=21$          & $n=41$          & $n=81$          \\
  \hline
  $m=32$ & $4.44 \pm 0.49$ & $1.59 \pm 0.15$ & $0.83 \pm 0.03$ & $0.65 \pm 0.03$ \\
  $m=16$ & $4.27 \pm 0.68$ & $1.82 \pm 0.09$ & $0.97 \pm 0.06$ & $0.83 \pm 0.03$ \\
  $m=8$  & $4.10 \pm 0.68$ & $2.24 \pm 0.21$ & $1.27 \pm 0.07$ & $1.14 \pm 0.03$ \\
  $m=4$  & $4.95 \pm 0.32$ & $2.77 \pm 0.19$ & $1.90 \pm 0.06$ & $1.78 \pm 0.04$ \\
  \hline
\end{tabular}
\end{center}
\caption{The error of the discrete Urysohn model in $\%$ depending on the number of rows, $m$, which corresponds to the discretisation of the integral operator in time, and the number of columns, $n$, which corresponds to the discretisation of the control function, of the Urysohn matrix. The averages and the $95\%$ confidence interval are calculated based on $9$ replications.}
\label{tab:errSto}
\end{table}

\subsubsection{System with noisy output}
\label{sec:resultsNoisyO}

The case that has been considered in section \ref{sec:resultsDiscr} is the ideal case when the output signal does not contain an observation error. However, in real systems, the output signal is often affected by the noise. In this case, the identification procedure should be robust and still identify the Urysohn matrix accurately.

To demonstrate the effect of the noise on the identification procedure, the same control function as in section \ref{sec:resultsDiscr} is considered. Moreover, the same identification and the same validation procedures are employed. However, after the output signal, which is subsequently used for the identification, is calculated, it is mixed with the noise. Thus, the Urysohn matrices are identified based on the following output:
\begin{equation*}
  y_i^\mathrm{S} = y_i + y_\mathrm{smax} \sigma w_i ,
\end{equation*}
where $y_i$ is the exact solution and $w_i \sim \mathcal{N}\left( 0, 1 \right)$ are normally distributed random numbers with zero mean and unit variance. The total simulation time is taken to be $t_\mathrm{max} = 4 \cdot 10^4$. The identification of the Urysohn matrices is performed with different values of stabilisation parameter $\alpha$. 

The results for different noise levels and different values of the stabilisation parameter are presented in table \ref{tab:errNoise}. The noise levels are $5\%$, $10\%$ and $20\%$ of the maximum static value of the solution, $y_\mathrm{smax}$. It can be seen that even for relatively large noise levels, the resulting Uryshon matrix can capture the system behaviour accurately and result in an error around $1\%$ across the validation dataset. For large noise levels, a relatively small value of $\alpha$ must be used to achieve high accuracy.

\begin{table}
\begin{center}
\begin{tabular}{ | r | c c c | } 
  \hline
                & $\alpha = 0.01$ & $\alpha = 0.05$ & $\alpha = 0.25$ \\
  \hline
  $\sigma=0.05$ & $0.44 \pm 0.03$ & $0.70 \pm 0.03$ & $1.50 \pm 0.11$ \\
  $\sigma=0.10$ & $0.67 \pm 0.05$ & $1.34 \pm 0.10$ & $2.99 \pm 0.24$ \\
  $\sigma=0.20$ & $1.18 \pm 0.06$ & $2.59 \pm 0.16$ & $5.95 \pm 0.38$ \\
  \hline
\end{tabular}
\end{center}
\caption{The error of the discrete Urysohn model in $\%$ depending on the noise level, $\sigma$, and the stabilisation parameter, $\alpha$, for the system with the noisy output. The averages and the $95\%$ confidence interval are calculated based on $9$ replications.}
\label{tab:errNoise}
\end{table}

\subsubsection{System with noisy input and output}
\label{sec:resultsNoisyIO}

The identification procedure of the discrete Urysohn operator can also be applied to systems with the noisy input and output. The same control function, the same identification and the same validation procedures as in section \ref{sec:resultsCont} are employed. However, before the coarse-sampled signals for the identification are calculated, the noise is added to both input and output:
\begin{align*}
  &x_i^\mathrm{S} = x_i + x_\mathrm{smax} \sigma v_i , \\
  &y_i^\mathrm{S} = y_i + y_\mathrm{smax} \sigma w_i ,
\end{align*}
where $x_i$ and $y_i$ are the exact input and output, respectively; $v_i \sim \mathcal{N}\left( 0, 1 \right)$ and $w_i \sim \mathcal{N}\left( 0, 1 \right)$ are normally distributed random numbers with zero mean and unit variance; $x_\mathrm{smax} = 1$ is the maximum value of the control. The discretisation of the Urysohn operator with $m=32$ and $n=81$ is used. The total simulation time is taken to be $t_\mathrm{max} = 4 \cdot 10^4$.

The results for different noise levels and different values of the stabilisation parameter are presented in table \ref{tab:errNoiseSto}. The noise levels are $5\%$, $10\%$ and $20\%$. As for the case of the system with the noisy output, the result indicates that the accuracy of the identified Urysohn model can easily be controlled by parameter $\alpha$ and the error decreases with the decrease of $\alpha$.

\begin{table}
\begin{center}
\begin{tabular}{ | r | c c c | } 
  \hline
                & $\alpha = 0.05$  & $\alpha = 0.20$  & $\alpha = 0.80$ \\
  \hline
  $\sigma=0.05$ & $0.81 \pm 0.02$  & $1.03 \pm 0.05$  & $2.33 \pm 0.15$ \\
  $\sigma=0.10$ & $1.44 \pm 0.10$  & $2.07 \pm 0.11$  & $4.72 \pm 0.27$ \\
  $\sigma=0.20$ & $3.89 \pm 0.16$  & $4.77 \pm 0.26$  & $9.43 \pm 0.44$ \\
  \hline
\end{tabular}
\end{center}
\caption{The error of the discrete Urysohn model in $\%$ depending on the noise level, $\sigma$, and the stabilisation parameter, $\alpha$, for the system with the noisy input and output. The averages and the $95\%$ confidence interval are calculated based on $9$ replications.}
\label{tab:errNoiseSto}
\end{table}

\subsection{Identification of real objects using experimental data}
\label{sec:grid}

The performance of the discrete Urysohn model and the proposed identification procedure has also been tested on real, experimentally-obtained data. Due to the scope of this paper, only the major results of this identification are briefly summarised in this section.

The identification of real objects has been conducted using publicly available datasets \cite{GroundVibration}, which were specifically collected and published for benchmarking identification algorithms. The most challenging task for the identification has been F-16 ground vibration test \cite{Noel2017}, for which the experimental data had been acquired on a full-scale F-16 aircraft. The tested object has non-linearities of clearance and friction type. Furthermore, the object has two inputs and multiple outputs.

The recorded data represent long series ($73$K to $116$K values), which is divided into two subsets --- the first subset has been recommended for the identification, the second subset has been recommended for the validation of the accuracy of the model on the unseen data. In contrast to the first example considered in this paper, this object has two inputs; therefore, a three-dimensional Urysohn kernel must be used to describe it. The identification procedure has been performed on several separate datasets. The number of time layers of the discrete Urysohn operators was taken between $160$ and $200$, while the number of quantisation levels for both inputs was taken from $80$ to $100$.

The validation of the obtained discrete Urysohn models has been performed on separate datasets, which were not involved in the identification process. The major result of the validation tests is that the observed error for different outputs has been between $1.2\%$ and $2.5\%$. For this comparison, the error is defined as the $L^2$-norm of the difference between the vector-columns of the recorded and the modelled outputs, which is divided by the difference between the maximum and the minimum values of the output and $\sqrt{N_\mathrm{p}}$, where $N_\mathrm{p}$ is the number of data points. Other details are omitted here and can be obtained directly from the C\# source code for this example, which is available as the supplementary information at \cite{PolarURL}.

\section{Conclusions}
\label{sec:concl}

The discrete Urysohn operator is a very efficient tool for modelling non-linear control systems due to its descriptive capabilities, low computational complexity, simplicity of identification and possibility of implementation in a few lines of code. The discrete model is an approximation of the continuous Urysohn model by certain quadrature rules. In this paper, the convergence has been demonstrated using a numerical example --- the modelling error decreased with the decrease of the discretisation step.

It has been shown that the model can be identified using an iterative algorithm based on the Kaczmarz method for solving linear systems of equations. The identification algorithm is simple and computationally inexpensive. Furthermore, the identification can be preformed using only observable data, i.e. the prescription of a specific input function is not required. The method improves the model as new input-output data points become available, without collecting long input and output sequences, and thus making the algorithm suitable for a real-time identification. The method does not require computational power for solving linear systems either. Therefore, the method is ideal for implementation in microelectronic systems with limited resources, such as battery management systems (BMS).

The noise of the input and the output data can be filtered out using the stabilisation parameter, which has been introduced in the identification procedure. Using a numerical example, it has been shown that by reducing the value of this parameter, the accuracy of the obtained Uryshon matrix increases, although the number of iterations that are required to obtain the accurate matrix also increases.

Since the Urysohn model is a general case of the Hammerstein and the linear models, the proposed identification method also covers all three nested models --- the Urysohn, the Hammerstein and the linear models. Usually, the suitability criteria for the choice of the model is built on the comparison of the computed and the measured outputs. In the case of the Uryshon model, the judgement can be made based on the values of the Urysohn matrix, from which it can be verified whether the system is describable by either the Hammerstein or the linear model\footnote{Discussed in appendix \ref{sec:relHammer}}.

This paper focused mainly on the simplest case --- the discrete operator with a quantised input. However, a generalisation of the model and the identification method for the case of a non-quantised input were also proposed. In this case, the number of parameters required to describe the system can be significantly reduced. 

The C\# implementation of the MISO model and the identification algorithm, as well as the demo programme for the computational example of section \ref{sec:grid} are available at \cite{PolarURL}. The minimalistic Matlab codes for the SISO model and its identification are provided in appendix \ref{sec:codes}.

\appendix
\titleformat{\section}[hang]{\Large\bfseries\raggedright\sffamily}{Appendix \thesection}{1em}{}

\section{Non-uniqueness of the Urysohn matrix}
\label{sec:rank}

\begin{proof}
The discrete Urysohn operator is always applied to $m$ successive elements of the input sequence and results in a single element of the output sequence. Therefore, the consideration of the single arbitrary input sequence and the corresponding output sequence of an Urysohn system is equivalent to the consideration of all possible inputs and outputs of the discrete Urysohn operator.

According to \eqref{eq:UryshonDiscr}, the input sequence of the discrete Urysohn operator consists of $m$ integer numbers, each of which can take a value from $1$ to $n$. The output of the discrete Urysohn operator is a single value. All possible outputs of the operator are denoted as $y_i^*$ and all possible inputs of the operator are denoted as sequences $K^i$. Elements of these sequences are denoted as $K^i_j$. Here $i \in \left\lbrace 1, 2, \ldots, N \right\rbrace$, where $N = n^m$.

Elements of $U$ and $y_i^*$ are rearranged into columns:
\begin{align}
  &Z = \begin{bmatrix}
    U_{11} & \ldots & U_{1n} & U_{21} & \ldots & U_{2n} & \ldots & U_{m1} & \ldots & U_{mn}
  \end{bmatrix}\transp ,
  \label{eq:UrysohnColApp} \\
  &Y = \begin{bmatrix}
    y_{1}^* & y_{2}^* & \ldots & y_{N}^*
  \end{bmatrix}\transp .
\end{align}
In this case, $Y$ is the product of matrix $M^m$ and column $Z$:
\begin{equation}
  M^m Z = Y ,
  \label{eq:UrysohnSLAE_App}
\end{equation}
where the elements of $M^m$ are given by
\begin{equation}
  M^m_{iq} = \begin{cases}
    1, & \text{if } q = n\left(j-1\right) + K^i_{m-j+1} \text{ where } j\in\left\lbrace 1, 2, \ldots, m \right\rbrace \\
    0, & \text{otherwise.}
  \end{cases}
  \label{eq:UrysohnColMatr_App}
\end{equation}
Such structure of matrix $M^m$ results directly from \eqref{eq:UryshonDiscr} and \eqref{eq:UrysohnColApp}. Matrix $M^m$ is formed by all possible input sequences and has $m n$ columns. Superscript $m$ in $M^m$ indicates the size of the input sequences. To prove the theorem, it must be proved that \textbf{(I)} rank of $M^m$ is $\left(m n - m + 1\right)$ and \textbf{(II)} for a selection of $\left(m n - m + 1\right)$ columns of $M^m$ to be linearly independent, it is necessary that at least $\left(n-1\right)$ columns are taken from each block of columns $\left(nj-n+1\right),\ldots,nj$, where $j\in\left\lbrace 1, 2, \ldots, m \right\rbrace$. From statement \textbf{(I)}, it follows that the Urysohn matrix is non-unique and has $\left(m-1\right)$ free parameters, as the elements of $U$ form the solution of \eqref{eq:UrysohnSLAE_App}. From statement \textbf{(II)}, it follows that each row of the Urysohn matrix cannot contain more than $1$ free parameter\footnote{Suppose row $b$ of the Urysohn matrix contains $2$ free parameters out of $\left(m-1\right)$, which implies that the corresponding columns of $M^m$ are linearly dependent on the other columns. Then, for a selection of $\left(m n - m + 1\right)$ linearly independent columns of $M^m$, the maximum of $\left(n-2\right)$ columns can be taken from block $\left(nb-n+1\right),\ldots,nb$. This contradicts statement \textbf{(II)}.}. The linear dependence of remaining $\left(m n - m + 1\right)$ elements of $U$ on $\left(m-1\right)$ parameters follows from linear system \eqref{eq:UrysohnSLAE_App}.

Statement \textbf{(I)} is proved by induction. If $m = 1$, the input sequences for the discrete Urysohn operator consist of a single integer number. Moreover, all possible input sequences are just numbers from $1$ to $n$. Thus, according to \eqref{eq:UrysohnColMatr_App}, $M^1$ is either the identity matrix of size $n$ or matrix, which is obtained from the identity matrix by the rearrangement of rows. Therefore, for $m = 1$, matrix $M^m$ has rank $n$.

Now it must be shown that if matrix $M^m$ has exactly $\left(m n - m + 1\right)$ independent rows for the Urysohn matrix of size $m \times n$, which was denoted above as $U$, then matrix $M^{m+1}$ has exactly $\left(m n - m + n\right)$ independent rows for the Urysohn matrix of size $\left(m+1\right) \times n$, which is denoted as $\bar{U}$.

Due to the inductive assumption, it can also be assumed that matrix $M^m$ that corresponds to all possible inputs of the operator with matrix $U$ is already assembled. It is useful to define $G^q$ to be a rectangular matrix where elements of column $q$ are equal to $1$, while all other elements are equal to $0$. Matrix $G^q$ has $n$ columns and the same number of rows as matrix $M^m$.

For the operator with matrix $\bar{U}$, the input sequences are longer by one number. Obviously, if $K^i$ are all possible inputs of the operator with matrix $U$, then 
\begin{equation*}
  \begin{bmatrix}
    s & K^i 
  \end{bmatrix} , \quad s\in\left\lbrace 1, 2, \ldots, n \right\rbrace
\end{equation*}
are all possible inputs of the operator with matrix $\bar{U}$. Therefore, according to \eqref{eq:UrysohnColMatr_App}, matrix $M^{m+1}$ that contains all possible inputs of the operator with matrix $\bar{U}$ is given by
\begin{equation}
  M^{m+1} = \begin{bmatrix}
    G^1 & M^m \\
    G^2 & M^m \\
    \vdots & \vdots \\
    G^n & M^m
  \end{bmatrix} .
  \label{eq:Minduct}
\end{equation}

The first block row (which contains matrices $G^1$ and $M^m$) has the number of independent rows of $\left(m n - m + 1\right)$ due to the inductive assumption. The second block row (which contains matrices $G^2$ and $M^m$) adds only one independent row, which is obvious from subtraction of the first block row from the second block row. The same is true for all remaining block rows. This results in
\begin{equation*}
  \left(m n - m + 1\right) + \left(n-1\right) = m n - m + n
\end{equation*}
independent rows of matrix $M^{m+1}$.

Thus, by induction, it has been proved that $M^m$ has exactly $\left(m n - m + 1\right)$ independent rows. By the fundamental theorem of linear algebra, this is also the rank of $M^m$. This concludes the proof of statement \textbf{(I)}.

To prove statement \textbf{(II)}, all columns of matrix $M^m$ are grouped into blocks $\left(nj-n+1\right),\ldots,nj$, where $j\in\left\lbrace 1, 2, \ldots, m \right\rbrace$. Suppose $\left(m n - m + 1\right)$ columns of $M^m$ are linearly independent and selected such that less than $\left(n-1\right)$ columns are taken from some block $j$. Since the total number of blocks is $m$, there are at least $2$ blocks, from which $n$ columns are taken. Without loss of generality, it can be assumed that these two blocks are the last two (otherwise, this can be achieved by a rearrangement of the corresponding blocks in $Z$). 

When the linear dependence of a set of columns of $M^m$ is considered, the order of rows in $M^m$ is irrelevant. In this case, $M^m$ can be formed incrementally according to equation \eqref{eq:Minduct}, i.e. $M^m$ formed from $M^{m-1}$, $M^{m-1}$ from $M^{m-2}$, etc. Thus, the block of the last two columns of $M^m$ is $M^2$ repeated vertically $n^{m-2}$ times. By statement \textbf{(I)}, the rank of $M^2$ is $\left(2n-1\right)$. This contradicts the assumed proposition, as the columns of the last two blocks (i.e. $2n$ columns) are supposed to be linearly independent. Therefore, the proposition is untrue, which concludes the proof of statement \textbf{(II)}.
\end{proof}

\begin{remark}
System \eqref{eq:UrysohnSLAE_App} is consistent, i.e. rank of 
\begin{equation*}
  \begin{bmatrix}
    M^m & Y 
  \end{bmatrix}
\end{equation*}
is also $\left(m n - m + 1\right)$, due to the conditional statement of the theorem --- the outputs are of the Urysohn system.
\end{remark}

\section{Convergence of the iterative method}
\label{sec:identConvProof}

\begin{proof}
As was shown in section \ref{sec:rank}, the Urysohn matrix contains $\left(m - 1\right)$ free parameters and the solution of the identification problem is non-unique. Within this proof, $U$ stands for any Urysohn matrix that exactly transforms the input sequence to the output sequence. The consequences of the existence of multiple solutions of the identification problem are unravelled closer to the end of the proof.

First of all, additional notation is introduced. Since the estimated Urysohn matrix changes each iteration, the iteration subscript is added to the estimated matrix and the model output is denoted as $\hat{y}_i$ to make the notation consistent with the description of the iterative algorithm,
\begin{equation*}
  \hat{y}_i = \sum_{j=1}^{m} U^i\left[ j, k_{i-j+1} \right] .
\end{equation*}
The actual recorded output is $y_i$, while $k_i$ is the input, which can be used instead of $x_i$ without loss of generality, according to \eqref{eq:UryshonDiscrContr}.

Elements of matrices $U$ and $U^i$ are rearranged into columns:
\begin{align*}
  &Z = \begin{bmatrix}
    U_{11} & \ldots & U_{1n} & U_{21} & \ldots & U_{2n} & \ldots & U_{m1} & \ldots & U_{mn}
  \end{bmatrix}\transp , \\
  &Z^i = \begin{bmatrix}
    U^i_{11} & \ldots & U^i_{1n} & U^i_{21} & \ldots & U^i_{2n} & \ldots & U^i_{m1} & \ldots & U^i_{mn}
  \end{bmatrix}\transp .
\end{align*}
The $L^2$-norm of $Z-Z^i$ is introduced and denoted as $e_i$,
\begin{equation}
  {e_i}^2 = \left( U_{11} - U^i_{11} \right)^2 + \left( U_{12} - U^i_{12} \right)^2 + \ldots + \left( U_{mn} - U^i_{mn} \right)^2 .
\end{equation}
At each iteration $i$ only a subset of $Z^i$ changes. Therefore, to simplify the notation,
\begin{align*}
  &a_j = U\left[ j, k_{i-j+1} \right] , \\
  &a^i_j = U^i\left[ j, k_{i-j+1} \right] ,
\end{align*}
are introduced. Elements $a^i_j$ are the only elements of $Z^i$, which are modified at iteration $i$.

By the iterative algorithm
\begin{equation}
  a^{i+1}_q =  a^i_q + \frac{\alpha}{m}\left( y_i - \hat{y}_i \right) = 
  a^i_q - \frac{\alpha}{m}\left( \sum_{j=1}^m a^i_j - \sum_{j=1}^m a_j \right) ,
  \quad q \in \left\lbrace 1,2,\ldots,m \right\rbrace .
\end{equation}
This leads to
\begin{multline}
  \sum_{q=1}^m \left( a^{i+1}_q - a_q \right)^2 = \\
  = \sum_{q=1}^m \left( \left( a^i_q - a_q \right)^2 - 2\frac{\alpha}{m}\left( \sum_{j=1}^m a^i_j -
  \sum_{j=1}^m a_j \right)\left( a^i_q - a_q \right) +
  \frac{\alpha^2}{m^2}\left( \sum_{j=1}^m a^i_j - \sum_{j=1}^m a_j \right)^2 \right) = \\
  = \sum_{q=1}^m \left( a^i_q - a_q \right)^2 +
  \frac{\alpha^2 - 2\alpha}{m}\left( \sum_{j=1}^m a^i_j - \sum_{j=1}^m a_j \right)^2 .
  \label{eq:errIter}
\end{multline}
Finally, since $a^i_j$ are the only elements of $Z^i$, which are modified at iteration $i$, equation \eqref{eq:errIter} leads to
\begin{equation}
  {e_{i+1}}^2 = {e_i}^2 - \frac{2\alpha - \alpha^2}{m}\left( \hat{y}_i - y_i \right)^2 .
  \label{eq:errDec}
\end{equation}
This means that for $\alpha \in \left(0, 2\right)$, if the model output is not equal to the exact output, the error $e_i$ necessarily decreases. Moreover, it can be seen that the fastest error decrease is achieved at $\alpha = 1$.

Up to now, the non-uniqueness of $U$ has not been used. Due to theorem \ref{th:ident}, all elements of matrix $U$ depend linearly on free $\left(m-1\right)$ parameters. Therefore, all possible solutions $Z$ of the identification problem form a flat $\left(m-1\right)$-size subspace of the $mn$-dimensional space. A numerical solution $Z^i$ is a particular point in the $mn$-dimensional space. Equation \eqref{eq:errDec} shows that the distance between $Z^i$ and all points of the flat subspace $Z$ decreases at each iteration, as long as $\hat{y}_i \neq y_i$. This could only mean that at each iteration, if $\hat{y}_i \neq y_i$ then $Z^i$ moves towards $Z$ in a direction, which is perpendicular to $Z$. This is schematically illustrated in figure \ref{fig:conv}.

\begin{figure}
  \begin{center}
    \includegraphics{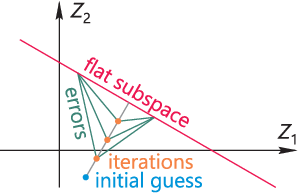}
  \end{center}
  \caption{A schematic illustration of the convergence of the iterative method. For explanatory purposes, $Z$ of size $2$ is taken. A flat $\left(m-1\right)$-size subspace is shown in red colour. An initial guess for $Z^i$ is shown in blue colour. The distances between $Z^i$ and two arbitrarily selected points of the flat subspace are shown in green colour. Since distances to all points of the flat subspace decrease at each iteration, $Z^i$ can move only in the perpendicular direction to the flat subspace.}
  \label{fig:conv}
\end{figure}

Due to the conditional statement of the theorem, the input sequence for the discrete Urysohn operator changes almost always. Thus, $\hat{y}_i \neq y_i$ almost always, since the exact input and output sequences are considered. Therefore, the iterative method converges to a solution, which has the minimum possible distance to the initial approximation, $U^m$, in the $mn$-dimensional space. Thus, when the starting point for the iterative algorithm is the all-zero matrix, the method converges to a unique solution, for which $Z$ has the minimum $L^2$-norm.
\end{proof}

\section{Describability by the discrete Urysohn operator}
\label{sec:critProof}

\begin{proof}
To prove the theorem, the notation is extended in the following way:
\begin{equation}
  y_i = y_i \left( X^i \right) , \quad
  X^i = \begin{bmatrix}
    x_{i-m+1} & \ldots & x_{i-1} & x_i
  \end{bmatrix} ,
\end{equation}
where $y_i$ is denoted as a function of $m$ quantised input values from $x_{i-m+1}$ to $x_i$. Additional vector-rows are introduced: $X^*$ is the vector-row of size $m$, in which all elements are equal to $x_\mathrm{min}$, and $X^{*p}$ is the vector-row of size $m$, in which all elements are equal to $x_\mathrm{min}$ except element $p$, which can take an arbitrary quantised value. Using this notation, part (b) implies the following:
\begin{equation}
  y_i \left( X^* \right) + y_i \left( X^{*p} + X^{*q} - x_\mathrm{min} \right) = y \left( X^{*p} \right) +
  y \left( X^{*q} \right) , \quad
  \forall p, q, \quad p \neq q .
\end{equation}

\textbf{Necessary condition.} Part (a) directly follows from the structure of the discrete model, equation \eqref{eq:UryshonDiscr}, since output $y_i$ is completely determined by input $x_i$ and preceding $\left(m-1\right)$ input values.

Part (b) follows from theorem \ref{th:ident}. Indeed, $\left(m - 1\right)$ elements of the Urysohn matrix, selected such that not more than one from each row is taken, can be prescribed arbitrary values. This means that, without loss of generality, elements $\left[2,1\right]$ to $\left[m,1\right]$ of the Urysohn matrix can be prescribed to be equal to element $\left[1,1\right]$. Furthermore, following the definition of $X^{*p}$ and definition \eqref{eq:UryshonDiscrContr} of the quantised input, the $p$-th element of $X^{*p}$ can be expressed as
\begin{equation*}
  x_\mathrm{min} + \left( k^{*p} - 1 \right) \frac{x_\mathrm{max} - x_\mathrm{min}}{n-1} ,
\end{equation*}
where $k^{*p}$ is an integer from $1$ to $n$. Analogously, an expression for the $q$-th element of $X^{*q}$ can be written. Without loss of generality $p<q$ can be taken. Finally, it can be seen that 
\begin{multline}
  y_i \left( X^* \right) + y_i \left( X^{*p} + X^{*q} - x_\mathrm{min} \right) = \sum_{j=1}^m U\left[ j, 1 \right] + 
  \sum_{j=1}^{m-p} U\left[ j, 1 \right] + \\ + 
  U\left[m-p+1,k^{*p}\right] + \sum_{j=m-p+2}^{m-q} U\left[ j, 1 \right] +
  U\left[m-q+1,k^{*q}\right] + \sum_{j=m-q+2}^m U\left[ j, 1 \right] = \\
  = \left( 2 m - 2 \right) U\left[ 1, 1 \right] + U\left[m-p+1,k^{*p}\right] + U\left[m-q+1,k^{*q}\right] = \\
  = \sum_{j=1}^{m-p} U\left[ j, 1 \right] + U\left[m-p+1,k^{*p}\right] + \sum_{j=m-p+2}^{m} U\left[ j, 1 \right] + \\
  + \sum_{j=1}^{m-q} U\left[ j, 1 \right] + U\left[m-q+1,k^{*q}\right] + \sum_{j=m-q+2}^{m} U\left[ j, 1 \right] =
  y \left( X^{*p} \right) + y \left( X^{*q} \right) .
\end{multline}

\textbf{Sufficient condition.} The sufficiency can be proved by first constructing an Urysohn matrix and then showing that any output will be described by an operator with such matrix. Given condition (a) and using the fact that the input is quantised, it follows that the system must be described by at most $mn$ parameters. Assume the following values for matrix $U$. The first column of matrix $U$ is assigned to be
\begin{equation*}
  U\left[1,1\right] = \ldots = U\left[m,1\right] = y_i \left( X^* \right) \frac{1}{m} .
\end{equation*}
Using expression for the $p$-th element of $X^{*p}$, which was introduced in the first part of the proof, all other elements are assigned to be
\begin{equation*}
  U\left[m-p+1,k^{*p}\right] = y_i \left( X^{*p} \right) - y_i \left( X^* \right) \frac{m-1}{m} , \quad
  k^{*p} \in \left\lbrace 2, \ldots, n \right\rbrace , \quad p \in \left\lbrace 1, \ldots, m \right\rbrace .
\end{equation*}
Now any input/output relationship of the system can be described by the constructed Urysohn operator. Indeed, consider output $y_i$ and input sequence $X^i$ given by
\begin{equation*}
  y_i = y_i \left( X^i \right) , \quad
  X^i = \begin{bmatrix}
    x_{i-m+1} & \ldots & x_{i-1} & x_i
  \end{bmatrix} .
\end{equation*}
Input sequence can be rewritten as
\begin{equation*}
  X^i = \begin{bmatrix}
    x_{i-m+1} \\ x_\mathrm{min} \\ x_\mathrm{min} \\ \vdots \\ x_\mathrm{min} \\ x_\mathrm{min}
  \end{bmatrix}\transp + 
  \begin{bmatrix}
    x_\mathrm{min} \\ x_{i-m+2} \\ x_\mathrm{min} \\ \vdots \\ x_\mathrm{min} \\ x_\mathrm{min}
  \end{bmatrix}\transp + \ldots +
  \begin{bmatrix}
    x_\mathrm{min} \\ x_\mathrm{min} \\ x_\mathrm{min} \\ \vdots \\ x_\mathrm{min} \\ x_i
  \end{bmatrix}\transp -
  \left(m-1\right) x_\mathrm{min} .
\end{equation*}
By applying condition (b) $m-1$ times, the following is obtained
\begin{multline*}
  y_i \left( X^i \right) = y_i \left(\begin{bmatrix}
    x_{i-m+1} \\ x_\mathrm{min} \\ \vdots \\ x_\mathrm{min}
  \end{bmatrix}\transp\right) + \ldots +
  y_i \left(\begin{bmatrix}
    x_\mathrm{min} \\ x_\mathrm{min} \\ \vdots \\ x_i
  \end{bmatrix}\transp\right) -
  \left(m-1\right) y_i \left( X^* \right) = \\
  = U\left[m,k_{i-m+1}\right] + \ldots + U\left[1,k_i\right] .
\end{multline*}
This shows that the input/output relationship is of the Urysohn system \eqref{eq:UryshonDiscr}.
\end{proof}

\section{Relation of the Urysohn model to other control systems' models}
\label{sec:relAll}

\subsection{Relation to the Hammerstein and the linear models}
\label{sec:relHammer}

One of the most well-known non-linear models of control systems is the Hammerstein model \cite{Hritonenko1996}:
\begin{equation}
  y\left( t \right) = \int_{0}^{T} h\left(s\right) H\left( x\left(t - s\right) \right) \dif s .
  \label{eq:HammerCont}
\end{equation}
In the case of block-model representation, the Hammerstein model is usually described by two blocks --- the first corresponding to a non-linear static part and the second corresponding to a linear dynamic part. By comparing equations \eqref{eq:UryshonCont} and \eqref{eq:HammerCont}, it can be seen that an Urysohn model, kernel $V\left( s, x \right)$ of which can be decomposed into a product of functions $h\left( s \right)$ and $H\left( x \right)$, is the Hammerstein model. When, in addition to that, function $H\left( x \right)$ is linear, the operator turns into the well-known convolution-type linear integral operator:
\begin{equation}
  y\left( t \right) = \int_{0}^{T} h\left(s\right) x\left(t - s\right) \dif s .
  \label{eq:linCont}
\end{equation}
For control systems described by the convolution-type linear integral operator, function $h\left(s\right)$ is the impulse response function.

It is useful to note that for a differentiable Urysohn kernel, it is possible to construct a linear approximation with respect to $x$ within a small variation of input $x$. Kernel smoothness is a usual property of many physical objects; therefore, in the case of a variation of the input within a small interval, the Urysohn model can be approximated by the linear model.

In the case of the discrete Urysohn model, it is easy to verify when the discrete Urysohn model becomes the discrete-time Hammerstein model or the discrete-time linear model. Matrix $U$ can always be expressed as a sum of matrices of the first rank, each of which is an outer product of two vectors. If such sum has only one significant term, i.e. all other terms can be neglected due to their order, matrix $U$ becomes the outer product of two vectors. For model \eqref{eq:UryshonDiscr}, $U\left[ j, k_{i-j+1} \right]$ becomes $\tilde{h}\left[ j \right] \tilde{H}\left[ k_{i-j+1} \right]$, where $\tilde{h}$ and $\tilde{H}$ are vector-columns, indices of which are shown in $\left[\cdot\right]$. Such model is the discrete-time quantised-input Hammerstein model, analogously to equation \eqref{eq:HammerCont}. If elements of $\tilde{H}$ linearly change with the index, the model becomes the discrete-time quantised-input linear model, analogously to equation \eqref{eq:linCont}.

\subsection{Relation to the Volterra series}
\label{sec:relVolt}

The continuous Urysohn operator is a particular case of the continuous-time Volterra series. Indeed, the general form of the continuous-time Volterra series is given by
\begin{equation}
  y\left(t\right) = h_0 + \sum_{q=1}^P \int_a^b \ldots \int_a^b h_q\left(\tau_1,\ldots,\tau_q\right) x\left(t-\tau_1\right) \ldots x\left(t-\tau_q\right) \dif \tau_1 \ldots \dif \tau_q .
  \label{eq:volter}
\end{equation}
If function $V\left(s,x\right)$ is smooth, $V$ in \eqref{eq:UryshonCont} can be Taylor-expanded with respect to the second variable, which leads to
\begin{equation}
  y\left(t\right) = \int_0^T \sum_{q=0}^\infty \frac{1}{q!}\left.\frac{\partial^q V\left(s,u\right)}{\partial u^q}\right|_{u=0} \left(x\left(t-s\right)\right)^q \dif s .
\end{equation}
This means that when
\begin{align}
  &h_q\left(\tau_1,\ldots,\tau_q\right) = \frac{1}{q!}\left.\frac{\partial^q V\left(\tau_1,u\right)}{\partial u^q}
    \right|_{u=0} \delta\left(\tau_1-\tau_2\right) \ldots \delta\left(\tau_1-\tau_q\right) , \quad q \geq 1 , \\
  &h_0 = \int_0^T V\left(s,0\right) \dif s , \quad
    a=0 , \quad b=T ,
\end{align}
and $P$ is infinity, the continuous-time Volterra model \eqref{eq:volter} becomes the continuous Urysohn model \eqref{eq:UryshonCont}.

The major limiting factor for efficient identification of the Volterra series model is its size. Usually, the model is limited to a relatively small number of terms (either in time domain or in frequency domain). Moreover, additional simplifications are often introduced (i.e. reduction of the number of the parameters to be identified), e.g. \cite{Schmidt2014}. For the exhaustive overview of the methods the reader is referred to \cite{Cheng2017}.

\subsection{Relation to the NARMAX model}
\label{sec:relNarm}

The discrete Urysohn operator is a particular case of the NARMAX model \cite{Chen1989}. Indeed, the general form of the NARMAX model is given by
\begin{equation}
  y_i = F\left( y_{i-1}, \ldots, y_{i-p}, x_{i-d}, \ldots, x_{i-d-m+1}, e_{i-1}, \ldots, e_{i-g} \right) + e_i ,
\end{equation}
where $x$, $y$, and $e$ are the input, the output and the error sequences, respectively; $F$ is a non-linear function. It is easy to verify that when $d=0$, $F$ does not depend on $y$ and $e$, and $F$ is given by
\begin{align}
  &F\left( x_{i}, \dots, x_{i-m+1} \right) = F_1\left( x_{i} \right) + \ldots + F_m\left( x_{i-m+1} \right) , \\
  &F_j\left(x_{i-j+1}\right) = U\left[ j, \frac{x_\mathrm{max} - n x_\mathrm{min} +
    \left(n-1\right)x_{i-j+1}}{x_\mathrm{max} - x_\mathrm{min}} \right] ,
\end{align}
the NARMAX model becomes the discrete Urysohn operator.

Although the discrete Urysohn model is only a particular case of the general NARMAX model, the latter is often simplified and polynomial expansions of the NARMAX model are usually used for identification/modelling purposes. Moreover, as in the case of the Volterra series, the set of the unknown parameters is often narrowed down (so-called ``structure detection''). Thus, identification methods may require intervention into the algorithms and expert knowledge of the underlying modelling system.

\section{Matlab codes}
\label{sec:codes}

The implementation of the discrete Urysohn model is relatively simple. The following Matlab function calculates the output of the Urysohn model based on the input:
\begin{lstlisting}[language=Matlab]
function [ y_ury ] = modelUrysohn( x, U, x_min, x_max )
%MODELURYSOHN Calculate the output of the Urysohn system based on the input
m = size(U,1);
n = size(U,2);
k = 1 + round( (n-1)*( x - x_min )/( x_max - x_min ) );
N = max(size(x));
y_ury = zeros(size(x));
for ii=m:N
    ctrl = k( ii:(-1):(ii-m+1) );
    ind = sub2ind( size(U), 1:m, ctrl );
    y_ury(ii) = sum( U(ind) );
end
end
\end{lstlisting}
Here \verb|x| is the input sequence; \verb|U| is the Uryshon matrix; \verb|x_min| and \verb|x_max| are the maximum and the minimum values of the input variable, respectively; and \verb|y_ury| is the calculated output sequence.

The iterative procedure for identifying the Urysohn matrix, which is suggested in section \ref{sec:identIter}, can be implemented in the following way:
\begin{lstlisting}[language=Matlab]
function [ U ] = identUrysohn( x, y, m, n, alpha, x_min, x_max )
%IDENTURYSOHN Identify the Urysohn matrix based on the input and the output
U = zeros(m,n);
k = 1 + round( (n-1)*( x - x_min )/( x_max - x_min ) );
N = max(size(x));
for ii=m:N
    y_real = y(ii);
    ctrl = k( ii:(-1):(ii-m+1) );
    ind = sub2ind( size(U), 1:m, ctrl );
    y_ury = sum( U(ind) );
    dy = y_real - y_ury;
    U(ind) = U(ind) + alpha*dy/m;
end
end
\end{lstlisting}
Here \verb|x| and \verb|y| are the input and the output sequences, respectively; \verb|m| and \verb|n| are the number of rows and columns of the Uryshon matrix, respectively; \verb|alpha| is the stabilisation parameter; \verb|x_min| and \verb|x_max| are the maximum and the minimum values of the input variable, respectively; and \verb|U| is the Uryshon matrix. For simplicity, the stopping criteria is not used in this implementation, and the iterative procedure loops until the end of the input sequence. Obviously, the stopping criteria can be added based on values of \verb|dy| for multiple consecutive iterations.

\begin{flushleft}

\bibliographystyle{unsrt}
\bibliography{refs}

\end{flushleft}

\end{document}